\documentclass[12pt]{article}
\usepackage{amsmath,graphicx,amssymb,amsthm}
\usepackage{rotating}
\allowdisplaybreaks
\def\cc{{\mathbb C}}

\def\U{\mathcal U}

\def\H{\mathcal H}

\def\L{\mathcal L}

\def\N{\mathcal N}

\def\Tr{{\rm Tr}}

\def\amslatex{$\mathcal{A}\kern-.1667em\lower.5ex\hbox{$\mathcal{M}$}\kern-.125em\mathcal{S}$-\LaTeX}

\def\tensor{\mathop{\bar\otimes}}
\newtheorem{set}{set}[section]
\newtheorem{Corollary}[set]{Corollary}

\newtheorem{Example}[set]{Example}
\newtheorem{Lemma}[set]{Lemma}

\newtheorem{Proposition}[set]{Proposition}
\newtheorem{Remark}[set]{Remark}
\newtheorem{Theorem}[set]{Theorem}
\newcommand{\define}{\mathrel{\hbox{$\equiv$\hskip -.90em \lower .47ex \hbox{$\leftharpoondown$}}}}
\newcommand{\enifed}{\mathrel{\hbox{$\equiv$\hskip -.90em \lower .47ex \hbox{$\rightharpoondown$}}}}

\numberwithin{equation}{section}
\begin{document}
\title{The Relative Weak Asymptotic Homomorphism\\ Property for Inclusions \\of Finite von Neumann Algebras}
\author{Junsheng Fang\and Mingchu Gao \and Roger R. Smith}
\date{}
 \maketitle
\begin{abstract}

A triple of finite von Neumann algebras $B\subseteq N\subseteq M$ is said to have the relative weak asymptotic homomorphism property if there exists a net of unitary operators $\{u_{\lambda}\}_{\lambda\in \Lambda}$ in $B$ such that $$\lim_{\lambda}\|{\mathbb{E}}_B(xu_{\lambda}y)-{\mathbb{E}}_B({\mathbb{E}}_N(x)u_{\lambda}{\mathbb{E}}_N(y))\|_2=0$$ for all $x,y\in M$.
 We prove that a triple of finite von Neumann algebras $B\subseteq N\subseteq M$
 has the relative weak asymptotic homomorphism property if and only if $N$ contains the  set of
 all $x\in M$ such that $Bx\subseteq \sum_{i=1}^n x_iB$ for a finite number of elements $x_1,\ldots,x_n$ in
 $M$. Such an $x$ is called a one sided quasi-normalizer of $B$, and the von Neumann algebra generated by all one sided quasi-normalizers of $B$ is called the one sided quasi-normalizer algebra of $B$.
 We characterize  one sided quasi-normalizer algebras  for  inclusions of group von Neumann algebras and use this to  show that  one sided quasi-normalizer algebras  and quasi-normalizer algebras  are not equal in general. We also give some applications to inclusions $L(H)\subseteq L(G)$ arising from containments of groups. For example, when $L(H)$ is a masa we determine the unitary normalizer algebra as the von Neumann algebra generated by the normalizers of $H$ in $G$.
\end{abstract}



\section{Introduction}
Let $M$ be a finite von Neumann algebra with a faithful normal
trace $\tau$, and let $B$ be a von Neumann subalgebra of $M$.  The
algebra $B$ has the \emph{weak asymptotic homomorphism property}
if there exists a net of unitary operators $\{u_{\lambda}\}_{\lambda\in \Lambda}$ in $B$
such that
\begin{equation}
\lim_{\lambda}\|{\mathbb{E}}_B(xu_{\lambda}y)-{\mathbb{E}}_B(x)u_{\lambda}{\mathbb{E}}_B(y)\|_2=0,\quad  x,y\in M.
\end{equation} 
This
property was  introduced by Robertson, Sinclair and
the third author~\cite{S-S,RSS,S-S2} for masas (maximal abelian subalgebras) of type ${\rm II}_1$ factors.
They showed that if $B$ has the weak asymptotic homomorphism
property, then $B$ is singular in $M$, and the purpose of introducing this property was to have an easily verifiable criterion for singularity. In~\cite{SSWW}, Sinclair,
 White, Wiggins and the third author proved that the converse is also true:
if $B$ is a singular masa  then $B$ has
the weak asymptotic homomorphism property (which is equivalent to the weakly mixing property of Jolissaint and Stalder \cite{J-S} for masas).  
We note that the results of \cite{SSWW} are formulated for $M$ a II$_1$ factor, but the proofs only require that $M$ be a finite von Neumann algebra with a faithful normal trace. Thus the equivalence of the weak asymptotic homomorphism property and singularity for masas is still valid at this greater level of generality. However, this equivalence breaks down beyond the masa case; in~\cite{G-W}, Grossman and Wiggins showed that if $B$ is a singular factor then $B$ does not necessarily have the weak asymptotic homomorphism property.

In~\cite{Ch}, Chifan introduced a generalized version as follows. A triple of von Neumann algebras $B\subseteq N\subseteq M$ is said to have the \emph{relative weak asymptotic homomorphism property} if there exists a net of unitary operators $\{u_{\lambda}\}_{\lambda\in \Lambda}$ in $B$ such that
\begin{equation}
 \lim_{\lambda}\|{\mathbb{E}}_B(xu_{\lambda}y)-{\mathbb{E}}_B({\mathbb{E}}_N(x)u_{\lambda}{\mathbb{E}}_N(y))\|_2=0,\quad   x,y\in M.
\end{equation}
Let $\N_M(B):=\{u\,\text{a unitary in}\, M:\,uBu^*=B\}$ denote the group of unitary normalizers of $B$ in $M$.
Chifan showed that if $B$ is a masa in a separable type ${\rm II}_1$ factor $M$, then
\begin{equation}
 B\subseteq\N_M(B)''\subseteq M
\end{equation}
has the relative weak asymptotic homomorphism property  (see~\cite{Muk} for a different proof).

A natural extension of Chifan's theorem is to consider a general
triple of finite von Neumann algebras $B\subseteq N\subseteq M$ and to ask for conditions which ensure that the relative weak asymptotic homomorphism property holds. Our main purpose in this paper is to provide a characterization of this property and to consider some subsequent applications. Our characterization is based on certain operators that are closely related to the quasi-normalizers introduced by Popa in \cite{Po1.5,Po2}.  Recall that he defined a {\em{quasi-normalizer}} for an inclusion $B\subseteq M$ to be an element $x\in M$ for which a finite set $\{x_1,\ldots,x_n\}\subseteq M$ can be found so that
\begin{equation}\label{eq1.1}
Bx\subseteq \sum_{i=1}^n x_i B,\ \ \ \ xB\subseteq \sum_{i=1}^n Bx_i,
\end{equation}
and we denote the set of quasi-normalizers by $q\mathcal{N}_M(B)$.  These are not quite the correct operators for our purposes, so we make a small adjustment by defining a {\em{one sided quasi-normalizer}} to be any element $x\in M$ satisfying only the first inclusion in \eqref{eq1.1}, and we denote the set of such elements by $q\mathcal{N}_M^{(1)}(B)$. 
In Section 3, we prove that a triple of finite von Neumann
algebras $B\subseteq N\subseteq M$ has the relative weak asymptotic
homomorphism property if and only if $N$ contains $q\mathcal{N}_M^{(1)}(B)$.   The von Neumann algebra generated by  $q\mathcal{N}_M^{(1)}(B)$ is called
the \emph{one sided quasi-normalizer algebra} of $B$ in  $M$, and is denoted by 
$W^*(q\mathcal{N}_M^{(1)}(B))$, a notation which reflects the fact that $q\mathcal{N}_M^{(1)}(B)$ is not necessarily self-adjoint. In the case that $B$ is a masa this characterization, combined with Chifan's theorem, gives that $W^*(q\mathcal{N}_M^{(1)}(B))=\mathcal{N}_M(B)''$. It had been shown earlier in \cite{P-S} that $q\mathcal{N}_M(B)''=\mathcal{N}_M(B)''$ when $M$ is a separable II$_1$ factor, and so $W^*(q\mathcal{N}_M^{(1)}(B))=q\mathcal{N}_M(B)''$ for masas, although these von Neumann algebras are different in general
(see Example \ref{ex5.3}). The advantage of the one sided quasi-normalizers is that they seem to be easier to calculate in specific examples, as we will see below. We note that one sided objects of this type play a significant role in understanding normalizers. For example, one sided unitary normalizers were important in \cite{SWW}, and a one sided version of groupoid normalizers was a key technical tool in \cite{FSWW}. 
 
After our characterization has been established in Section 3, we devote Section 4  to the study of $W^*(q\mathcal{N}_M^{(1)}(B))$. Here we show, among other results, that
 the one sided quasi-normalizer algebra and the
quasi-normalizer algebra of an atomic  von
Neumann subalgebra $B$ of a finite von Neumann algebra $M$ are equal to $M$.

In Section 5, we apply these results to inclusions of von Neumann algebras arising from inclusions $H\subseteq G$ of discrete groups.
We characterize  one sided quasi-normalizer algebras of such inclusions  in terms of properties of the groups, and also  show that  one sided quasi-normalizer algebras and
quasi-normalizer algebras are not equal in general. Making use of the one sided quasi-normalizers, we are able to study unitary normalizers and show, for example, that when $L(H)$ is a masa in $L(G)$, its unitary normalizer algebra is the von Neumann algebra of the group of normalizers of $H$ in $G$. This leads to new characterizations of when $L(H)$ is either singular or Cartan. In section 6, we summarize the relationships between various types of normalizer algebras, and we show that
   one sided quasi-normalizer algebras
have some special properties when compared to the other types. For example, we establish a tensor product formula in Proposition \ref{pro6.1} which parallels similar formulas for groupoid normalizers and intertwiners proved in \cite{FSWW}.

Sufficient background material on finite von Neumann algebras for this paper may be found in \cite{S-S2}.

{\bf Acknowledgement} This work originated during the Workshop on
Analysis and Probability held at Texas A$\&$M University in the
summer of 2009. It is our pleasure to thank the organizers of the
workshop and the NSF for the financial support to the workshop.  The authors also thank Ionut Chifan,  Kunal Mukherjee and Stuart White for valuable comments on this paper.

\section{Preliminaries}\label{sec2}

Throughout this paper, $M$ is a finite von Neumann algebra with a
given faithful normal trace $\tau$. We use $L^2(M)=L^2(M,\tau)$ to
denote the Hilbert space obtained by the GNS-construction of $M$
with respect to $\tau$.
  The image of $1\in M$ via the GNS construction is denoted by $\xi$ and the image of $x\in M$ is denoted by $x\xi$. Throughout this paper, we will reserve the letter $\xi$ for this purpose. The trace norm of $x\in M$ is defined by $\|x\|_{2}=\|x\|_{2,\tau}=\tau(x^*x)^{1/2}$.
The letter $J$ is reserved for the isometric conjugate linear operator on $L^2(M)$ defined on $M\xi$ by $J(x\xi)=x^*\xi$ and extended by continuity to $L^2(M)$ from the dense subspace $M\xi$.

Let $B\subseteq M$ be an inclusion of finite von Neumann algebras. Then there exists a unique
faithful normal conditional expectation ${\mathbb{E}}_B$ from $M$ onto $B$
preserving $\tau$. Let $e_B$ be the projection of $L^2(M)$ onto
$L^2(B)$. For $x\in M$, we have $e_B(x\xi)={\mathbb{E}}_B(x)\xi$.  The von Neumann algebra $\langle M,e_B\rangle$
generated by $M$ and $e_B$ is called \emph{the basic construction}, which plays a crucial role in the study of von Neumann
subalgebras of finite von Neumann algebras. The basic construction has many remarkable properties (see~\cite{Jones,P-P,S-S2}). In particular, there exists a unique
faithful tracial weight ${\rm Tr}$ on $\langle M,e_B\rangle$ such
that
\begin{equation}\label{E:Tr and tau}
 {\rm Tr}(xe_By)=\tau(xy),\quad  x,y\in M.
\end{equation}
Furthermore, we can choose a net of vectors $\{\xi_i\}_{i\in I}$ from $L^2(M)$  such that
\begin{equation}\label{E:Tr}
 \Tr(t)=\sum_{i\in I}\langle t\xi_i,\xi_i\rangle_{2,\tau},\quad t\in \langle M,e_B\rangle^+
\end{equation}
(see  \cite[Lemma 4.3.4, Theorem 4.3.11]{S-S2}).  An examination of the proof of \cite[Lemma 4.3.4]{S-S2} shows that we may construct the index set $I$ to have a minimal element $i=1$ and we may take $\xi_1$ to be $\xi$. Letting $t=e_B$ in~\eqref{E:Tr}, we have $e_B\xi_i=0$ for all $i\ne 1$, $i\in I$.

There is a well defined map $\Psi: Me_BM\rightarrow M$, given by
\begin{equation}
 \Psi(xe_By)=xy,\quad x,y\in M,
\end{equation}
and called the {\emph{pull down map}}.  It was shown in~\cite{P-P} (also see~\cite{S-S2}) that the pull down map can be extended to a contraction from $L^1(\langle M,e_B\rangle, {\rm Tr})$ to $L^1(M,\tau)$, which is just the predual of the embedding $M\hookrightarrow \langle M,e_B\rangle$.

Let $w\in \langle M,e_B\rangle$, and let $\eta=w(\xi)\in L^2(M)$. Then
\begin{equation}\label{E:def of Lxi}
 L_\eta(x\xi)=Jx^*J(\eta),\quad x\in M,
\end{equation}
is a densely defined operator affilated with $M$. We may identify $L_\eta$ with $\eta$ in a canonical way  so that $\|L_\eta\|_{2,\tau}=\tau(L_\eta^*L_\eta)^{1/2}=\|\eta\|_{2,\tau}$ is well defined (see~\cite{Ne}).    Note that $we_B=L_\eta e_B$. Indeed, for $x\in M$, we have
\begin{align}
 L_\eta e_B(x\xi)&=L_\eta({\mathbb{E}}_B(x)\xi)=J{\mathbb{E}}_B(x^*)J(\eta)
=J{\mathbb{E}}_B(x^*)Jw(\xi)\notag\\
&=wJ{\mathbb{E}}_B(x^*)J(\xi)=we_B(x\xi).
\end{align}

For $z\in \langle M,e_B\rangle$, define $\|z\|_{2,\Tr}=\Tr(z^*z)^{1/2}$.  The following lemmas are well known to experts.  For the  reader's convenience, we include the proofs.

\begin{Lemma}\label{L:pre 1}
 Suppose that $w\in \langle M,e_B\rangle$ and $\eta=w(\xi)\in L^2(M)$. Then
$\|we_B\|_{2,\Tr}=\|\eta\|_{2,\tau}$.
\end{Lemma}
\begin{proof}
 The equalities 
\begin{align}
 \|we_B\|_{2,\Tr}^2&=\Tr(e_Bw^*we_B)=\sum_{i\in I} \langle we_B\xi_i,we_B\xi_i\rangle_{2,\tau}\notag\\
&=\langle w(\xi),w(\xi)\rangle_{2,\tau}=\langle \eta,\eta\rangle_{2,\tau}=\|\eta\|_{2,\tau}^2
\end{align}
follow from \eqref{E:Tr} and the fact that $e_B\xi_i=0$ for $i\ne 1$.
\end{proof}

The following lemma plays a key role in the proof of  Lemma~\ref{L:Popa}.

\begin{Lemma}\label{L:pre 2}
 Suppose that $w\in \langle M,e_B\rangle$ and $\eta=w(\xi)\in L^2(M)$. Then
\begin{equation}
 \Psi(we_Bw^*)=L_\eta L_\eta^*,
\end{equation}
where $\Psi$ is the pull down map and $L_\eta$ is the operator defined by~\eqref{E:def of Lxi}.
\end{Lemma}
\begin{proof}
 Since $\eta \in L^2(M)$, there exists a sequence  $\{x_n\}_{n=1}^{\infty}$ in $M$ such that \begin{equation}
\lim_{n\rightarrow\infty}\|L_\eta-x_n\|_{2,\tau}=\lim_{n\rightarrow\infty}\|\eta-x_n\xi\|_{2,\tau}=0.
\end{equation} 
Therefore,
\begin{align}
\|L_\eta L_\eta^*-x_nx_n^*\|_{1,\tau}&\leq \|(L_\eta-x_n)L_\eta^*\|_{1,\tau}+\|x_n(L_\eta^*-x_n^*)\|_{1,\tau}\notag\\
&\leq
\|L_\eta-x_n\|_{2,\tau}\|L_\eta\|_{2,\tau}+\|x_n\|_{2,\tau}\|L_\eta-x_n\|_{2,\tau}\rightarrow 0.
\end{align}
By Lemma~\ref{L:pre 1},
\begin{equation}
\lim_{n\rightarrow\infty} \|we_B-x_ne_B\|_{2,\Tr}=\lim_{n\rightarrow\infty}\|\eta-x_n\xi\|_{2,\tau}=0,
\end{equation}
and so
$\lim_{n\rightarrow\infty}\|we_Bw^*-x_ne_Bx_n^*\|_{1,\Tr}$ $=0$. Noting that $\Psi(x_ne_Bx_n^*)=x_nx_n^*$, the equality $\Psi(we_Bw^*)$ $=$ $L_\eta L_\eta^*$ follows since  $\Psi$ is a continuous contraction from $L^1(\langle M,e_B\rangle, {\rm Tr})$ to $L^1(M,\tau)$.
\end{proof}

We have included here only the facts about the basic construction that we will need subsequently. 
Much more detailed coverage can be found in \cite{Chr, Jones, P-P, S-S2}.

\section{Main result}\label{sec3}

This section  is devoted to the main result of the paper, Theorem \ref{T:main theorem}. We state it immediately, but defer the proof until Lemmas \ref{L:chifan}--\ref{L:Popa} have been established. We have included part (iii) for emphasis, but it is of course just a notational restatement of (ii).

\begin{Theorem}\label{T:main theorem} The following conditions are equivalent for finite von Neumann algebras:
\begin{itemize}
\item[{\rm{(i)}}] The triple $B\subseteq N\subseteq M$ has the relative weak asymptotic homomorphism property;
\item[{\rm{(ii)}}] If $x\in M$ satisfies $Bx\subseteq \sum_{i=1}^n x_iB$ for a finite number of elements $x_1,\ldots,x_n$ in $M$, then $x\in N$;
\item[{\rm{(iii)}}] $q\mathcal{N}_M^{(1)}(B)\subseteq N$.
\end{itemize}
\end{Theorem}

To prove this theorem, we will need several lemmas. The following  is essentially  \cite[Lemma 2.5]{Ch} (see also  \cite[Corollary 2.3]{Po3}), and  so we omit the proof.
\begin{Lemma}\label{L:chifan}
Suppose that a triple of finite von Neumann algebras $B\subseteq N\subseteq M$ does not have the relative weak asymptotic homomorphism property. Then there exists a nonzero projection $p\in B'\cap \langle M,e_B\rangle$ such that $0<\Tr(p)<\infty$ and $p\leq 1-e_N$.
\end{Lemma}

\begin{Lemma}\label{L:approximate projections} Let $B\subseteq N\subseteq M$ be a triple of finite von Neumann algebras, let $p\in \langle M,e_B \rangle$ be a finite projection satisfying $p\leq 1-e_N$, and let
$\varepsilon>0$. Then there exists a finite number of elements
$x_1,\ldots,x_n\in M$ such that ${\mathbb{E}}_N(x_i)=0$ for $1\leq i\leq n$, and
\begin{equation}\label{eq3.z}
\left\|p-\sum_{i=1}^n x_ie_B x_i^*\right\|_{2,\Tr}<\varepsilon.
\end{equation}
\end{Lemma}
\begin{proof}
 By  \cite[Lemma 1.8]{Po1}, there are elements $y_1,\ldots, y_n\in M$  such that
\begin{equation}\label{E:2}
 \left\|p-\sum_{i=1}^n y_ie_B y_i^*\right\|_{2,\Tr}<\varepsilon/3.
\end{equation}
For $1\leq i\leq n$, let $x_i=y_i-{\mathbb{E}}_N(y_i)$, and note that ${\mathbb{E}}_N(x_i)=0$. Since $pe_N=0$, it follows from
\eqref{E:2} that
\begin{equation}\label{eq3.zz}
 \left\|\sum_{i=1}^n {\mathbb{E}}_N(y_i)e_By_i^*\right\|_{2,\Tr}=\left\|e_N\left(p-\sum_{i=1}^n y_ie_B y_i^*\right)\right\|_{2,\Tr}<\varepsilon/3.
\end{equation}
Also, the identity
\begin{align}
(1-e_N)\left(p-\sum_{i=1}^n y_ie_B y_i^*\right)e_N&=-\sum_{i=1}^n(1-e_N)y_ie_Ne_Be_Ny_i^*e_N\notag\\
&=-\sum_{i=1}^n(y_i-{\mathbb{E}}_N(y_i))e_B{\mathbb{E}}_N(y_i^*)\notag\\
&=-\sum_{i=1}^nx_ie_B{\mathbb{E}}_N(y_i^*)
\end{align}
shows that
\begin{equation}\label{eq3.zzz}
\left\|\sum_{i=1}^n x_ie_B{\mathbb{E}}_N(y_i^*)\right\|_{2,\Tr}=\left\|(1-e_N)\left(p-\sum_{i=1}^n y_ie_B y_i^*\right)e_N\right\|_{2,\Tr}<\varepsilon/3,
\end{equation}
from \eqref{E:2}.
Using the expansion
\begin{align}
y_ie_By_i^*&=(x_i+{\mathbb{E}}_N(y_i))e_By_i^*
=x_ie_By_i^*+{\mathbb{E}}_N(y_i)e_By_i^*\notag\\
&=x_ie_Bx_i^*+x_ie_B{\mathbb{E}}_N(y_i^*)+{\mathbb{E}}_N(y_i)e_By_i^*
\end{align}
and the inequalities \eqref{E:2}, \eqref{eq3.zz}, and \eqref{eq3.zzz}, we see that
\begin{align}
 \left\|p-\sum_{i=1}^n x_ie_B x_i^*\right\|_{2,\Tr}\leq &\left\|p-\sum_{i=1}^n y_ie_B y_i^*
\right\|_{2,\Tr}
 +\left\|\sum_{i=1}^n {\mathbb{E}}_N(y_i)e_By_i^*\right\|_{2,\Tr}\notag\\
 &+\left\|\sum_{i=1}^n x_ie_B\mathbb{E}_N( y_i^*)\right\|_{2,\Tr}<\varepsilon,
\end{align}
proving the result.
\end{proof}

If $p\in B'\cap \langle M,e_B\rangle$, then $\H=pL^2(M)$ is a $B$-bimodule. Conversely, if $\H\subseteq L^2(M)$ is a  $B$-bimodule and $p$ is the orthogonal projection of $L^2(M)$ onto $\H$, then $p\in B'\cap \langle M,e_B\rangle$. In the following we recall some basic facts about $B$-bimodules.

Suppose that a Hilbert subspace $\H\subseteq L^2(M)$ is a right $B$-module. Let $\L_B(L^2(B),\H)$  be the set of bounded right
$B$-modular operators from $L^2(B)$ into $\H$. For instance, if
$\H=L^2(B)$, then $\L_B(L^2(B),L^2(B))$ consists of operators
induced by the left action of $B$ on $L^2(B)$.

Let $B$ be a finite von Neumann algebra with a faithful normal trace $\tau$. Suppose that $B$ acts on the right on a Hilbert space $\H$. Then the dimension of $\H$ over $B$ is defined as
\begin{equation}
 {\rm dim}_B(\H)=\Tr(1),
\end{equation}
where $\Tr$ is the unique tracial weight on $B'$ satisfying the following condition
\begin{equation}
 \Tr(xx^*)=\tau(x^*x),\quad   x\in \L_B(L^2(B),\H).
\end{equation}
 We say that $\H$ is a \emph{finite right $B$-module} if $\Tr(1)<\infty$. For details of finite right $B$-modules,
 we refer to  \cite[Appendix A]{Va}.

Suppose that $\H\subseteq L^2(M)$ is a right $B$-module. Then $\H$ is called a
\emph{finitely generated right $B$-module} if there exists a finite set of 
elements $\{\eta_1,\ldots,\eta_n\}\subseteq \H$ such that $\H$ is the closure
of $\sum_{i=1}^n\eta_i B$.  A set $\{\eta_i\}_{i=1}^n$ of elements
in $\H$ is called an {\em{orthonormal basis}} of  $\H$ if
${\mathbb{E}}_B(\eta_i^*\eta_j)=\delta_{ij} p_i\in B$, where each $p_i$ is a projection and, for
every $\eta\in \H$, we have
\begin{equation}\label{eq3.ab}
\eta=\sum_{i=1}^n \eta_i {\mathbb{E}}_B(\eta_i^*\eta).
\end{equation}
Note that, by putting $\eta=\eta_j$ into \eqref{eq3.ab}, we have $\eta_j=\eta_jp_j$, $1\leq j\leq n$.
It might appear that the vectors on the right hand side of \eqref{eq3.ab} are not in $\H$,
since $\eta_i\in \H\subseteq L^2(M)$ and ${\mathbb{E}}_B(\eta_i^*\eta)\in L^1(B)$,
but the construction of the orthonormal basis ensures that they do lie in this Hilbert space.

 Let $p_\H$ be the orthogonal projection of $L^2(M)$ onto $\H$. Following \cite[Lemma 1.4.2]{Po2}, we have $p_\H=\sum_{i=1}^n L_{\eta_i}e_B L_{\eta_i}^*$, where $L_\eta$ is defined as in~\eqref{E:def of Lxi}. Let $w_i=L_{\eta_i}e_B$, a bounded operator since $w_iw_i^*\leq p_{\H}$. For each $x\in M$ and $b\in B$,
\begin{align}
 w_i JbJ (x\xi)&=w_i(xb^*\xi)=L_{\eta_i}e_B(xb^*\xi)=L_{\eta_i}({\mathbb{E}}_B(xb^*)\xi)\notag\\
&=L_{\eta_i}({\mathbb{E}}_B(x)b^*\xi)=Jb{\mathbb{E}}_B(x^*)J(\eta_i)\notag\\
&=JbJJ{\mathbb{E}}_B(x^*)Jw_i(\xi)=JbJw_iJ{\mathbb{E}}_B(x^*)J(\xi)=J bJw_i(x\xi).
\end{align}
 Thus $w_iJbJ=JbJw_i$, which implies that $w_i\in \langle M,e_B\rangle$. Summarizing the above arguments, we have shown that
\begin{equation}\label{E:unbounded}
 p_\H=\sum_{i=1}^n w_ie_Bw_i^*,
\end{equation}
where $w_i=L_{\eta_i}e_B\in\langle M,e_B\rangle$. We note that every finitely generated right $B$-module has an orthonormal basis, \cite[1.4.1]{Po2}.

The following lemma is proved by Vaes in~\cite[Lemma A.1]{Va} (see also \cite[Lemma 1.4.2]{Po2}). It is designed to circumvent the difficulty that finite right $B$-modules might not be finitely generated.

\begin{Lemma}\label{L:vaes}
 Suppose that $\H$ is a finite right $B$-module. Then there exists a sequence of projections 
$\{z_n\}_{n=1}^{\infty}$ in $Z(B)=B'\cap B$ and a sequence of integers $\{k_n\}_{n=1}^{\infty}$ such that $\lim_{n\rightarrow\infty}z_n=1$ in the strong operator topology and
 $\H z_n$ is unitarily equivalent to a  left $p_n{\mathbb{M}}_{k_n}(B)p_n$ right $B$-module  $p_n(L^2(B)^{(k_n)})$ for each $n$, where $p_n$ is a projection in ${\mathbb{M}}_{k_n}(B)$. In particular, $\H z_n$ is a finitely generated right $B$-module.
\end{Lemma}

The following lemma is motivated by  \cite[Lemma 1.4.2]{Po2}.
\begin{Lemma}\label{L:Popa}
Suppose that $\H\subseteq L^2(M)$ is a  $B$-bimodule, and that $\H$ is a finitely generated right $B$-module with an orthonormal basis of length $k$.   Let $p_\H$ be the orthogonal projection of $L^2(M)$ onto $\H$. Then there exists a sequence of   projections $z_n$ in $B'\cap M$  such that $\lim_{n\rightarrow\infty}z_n=1$ in the strong operator topology and for each $n$,
\begin{equation}
 z_np_\H z_n(x\xi)=\sum_{i=1}^kx_{n,i}{\mathbb{E}}_B(x_{n,i}^*x)\xi,\quad   x\in M,
\end{equation}
for a finite number of elements $x_{n,1},\ldots,x_{n,k}\in M$.
\end{Lemma}
\begin{proof}
Let $\{\eta_i\}_{i=1}^k\subseteq \H\subseteq L^2(M,\tau)$ be an orthonormal basis for $\H$, in which case $\H=\oplus_{i=1}^k [\eta_iB]$.  By~\eqref{E:unbounded},
 $p_\H=\sum_{i=1}^kw_ie_Bw_i^*\in B'\cap \langle M,e_B\rangle$, where $w_i=L_{\eta_i}e_B\in \langle M,e_B\rangle$.
For $b\in B$, we have
\begin{equation}
 b\sum_{i=1}^k w_ie_Bw_i^*=\sum_{i=1}^kw_ie_Bw_i^*b.
\end{equation}
Applying the pull down map to both sides and noting that $w_i(\xi)=\eta_i$,  we obtain
\begin{equation}
b\left(\sum_{i=1}^kL_{\eta_i}L_{\eta_i}^*\right)=\left(\sum_{i=1}^kL_{\eta_i}L_{\eta_i}^*\right)b
\end{equation}
by Lemma~\ref{L:pre 2}. Since $\sum_{i=1}^kL_{\eta_i}L_{\eta_i}^*$ is an operator affiliated with $M$,
$q\in B'\cap M$ for all spectral projections $q$ of $\sum_{i=1}^kL_{\eta_i}L_{\eta_i}^*$.
Therefore,  there exists  a sequence of projections $z_n\in B'\cap M$ such that $\lim_{n\rightarrow \infty}z_n=1$ in the strong operator topology and $\sum_{i=1}^kz_nL_{\eta_i}L_{\eta_i}^*z_n$ is a bounded operator for each $n$. Let $x_{n,i}=z_nL_{\eta_i}$, $1\leq i\leq k$. Then $x_{n,i}\in M$ and
\begin{equation}
 z_np_{\H}z_n(x\xi)=\sum_{i=1}^kx_{n,i}{\mathbb{E}}_B(x_{n,i}^*x)\xi,\quad   x\in M,
\end{equation}
as required.
\end{proof}

This completes the preparations for the proof of our main result, which we now give. We will establish only the equivalence of (i) and (ii) since, as already noted, (iii) is just a restatement of (ii).
\begin{proof}[Proof of Theorem~\ref{T:main theorem}.]
(ii)$\Rightarrow$(i).\quad To derive a contradiction, suppose that (ii) holds but that
 the triple $B\subseteq N\subseteq M$ does not have the relative weak asymptotic homomorphism property. By Lemma~\ref{L:chifan}, there exists a
 projection $p\in B'\cap \langle M,e_B\rangle$ such that $0<\Tr(p)<\infty$ and $p\leq 1-e_N$.

Let $\H=pL^2(M)$. Then $\H$ is a $B$-bimodule and a finite right $B$-module.  By Lemma~\ref{L:vaes}, we may assume that $\H$ is a finitely generated right $B$-module. By Lemma~\ref{L:Popa}, there exists a sequence of   projections $z_n$ in $B'\cap M$  such that $\lim_{n\rightarrow\infty}z_n=1$ in the strong operator topology and for each $n$,
\begin{equation}
 z_npz_n(x\xi)=\sum_{i=1}^kx_{n,i}{\mathbb{E}}_B(x_{n,i}^*x)\xi \in M\xi,\quad   x\in M,
\end{equation}
for a finite number of elements $x_{n,1},\ldots,x_{n,k}\in M$.  Note that $z_npz_n\in B'\cap \langle M,e_B\rangle$.
Thus, for every $x\in M$,
\begin{equation}
 B\left({z_npz_n(x\xi)}\right)=(z_npz_n)(Bx\xi)\subseteq\sum_{i=1}^kx_{n,i}B\xi.
\end{equation}
Thus $z_npz_n(x\xi)\in N\xi\subseteq L^2(N)$ by the assumption (ii) of Theorem~\ref{T:main theorem}. Hence, for each $\eta\in L^2(M)$, $z_npz_n(\eta)\in L^2(N)$.
 Since $\lim_{n\rightarrow\infty}z_n=1$ in the strong operator topology, 
\begin{equation}
p(\eta)=\lim_{n\rightarrow\infty} z_npz_n(\eta)\in L^2(N),
\end{equation} 
and so  $p\leq e_N$. On the other hand, $p\leq 1-e_N$ and we arrive at the contradiction $p=0$. 

\medskip

\noindent (i)$\Rightarrow$(ii).\quad Suppose that $x\in M$ satisfies $Bx\subseteq \sum_{i=1}^nx_iB$ for a finite number of elements $x_1,\ldots,x_n$ in $M$, and let $\H$ be the closure of $BxB\xi$ in $L^2(M)$. Then $\H$ is a $B$-bimodule and $\H\subseteq L^2(\sum_{i=1}^nx_iB)$. Thus $\H$ is a finite right $B$-module. Let $p$ be the projection of $L^2(M)$ onto $\H$. Then $p\in B'\cap \langle M,e_B\rangle$ and $0<\Tr(p)<\infty$. We need only prove that $p\leq e_N$ since if this is the case, then $x\xi=p(x\xi)=e_N(x\xi)\in L^2(N)$, implying that $x\in N$.

Suppose that $e_Npe_N=p$ is not true. Then $(1-e_N)p\neq 0$. Replacing $p$ by a nonzero  spectral projection of $(1-e_N)p(1-e_N)$ corresponding to some interval $[c,1]$ with $c>0$, we may assume that $p$ is a nonzero subprojection of $1-e_N$.

Let $\varepsilon>0$. By Lemma~\ref{L:approximate projections}, there exists a finite set of elements
$\{x_1,\ldots,x_n\}\subseteq M$ such that ${\mathbb{E}}_N(x_i)=0$ and
\begin{equation}\label{eq3.1}
\|p-\sum_{i=1}^n x_ie_B x_i^*\|_{2,\Tr}<\varepsilon/2.
\end{equation}
Let $p_0=\sum_{i=1}^n x_ie_B x_i^*$. Since $p\in B'\cap \langle M,e_B\rangle$, $upu^*=p$ for all unitary operators $u\in B$. Thus
\begin{equation}
 \|up_0u^*-p_0\|_{2,\Tr}\leq \|u(p_0-p)u^*\|_{2,\Tr}+\|p_0-p\|_{2,\Tr}<\varepsilon,\quad  u\in \U(B).
\end{equation}
Therefore,
\begin{align}
 2\|p_0\|_{2,\Tr}^2&=\|up_0u^*-p_0\|_{2,\Tr}^2+2\Tr(up_0u^*p_0)\notag\\
&=\|up_0u^*-p_0\|_{2,\Tr}^2+2\sum_{1\leq i,j\leq n} \Tr(ux_ie_Bx_i^*u^*x_je_Bx_j^*)\notag\\
&\leq \varepsilon^2+2\sum_{1\leq i,j\leq n}\tau({\mathbb{E}}_B(x_i^*u^*x_j)x_j^*ux_i)\notag\\
&\leq \varepsilon^2+2\sum_{1\leq i,j\leq n}\|{\mathbb{E}}_B(x_j^*ux_i)\|_{2,\tau}^2
\end{align}
for all unitary operators $u$ in $B$.  By the assumption of (i), there exists a sequence of unitary operators $\{u_k\}_{k=1}^{\infty}$ in $B$ such that  $\sum_{1\leq i,j\leq n}\|{\mathbb{E}}_B(x_j^*u_kx_i)\|_{2,\tau}^2\rightarrow 0$ when $k\rightarrow \infty$. Hence, $\|p_0\|_{2,\Tr}< \varepsilon$.  Since $\varepsilon>0$ is arbitrary, it follows from \eqref{eq3.1} that $p=0$, giving a contradiction and completing the proof.
\end{proof}

\section{One sided quasi-normalizer algebras}\label{sec4}

Recall that an element $x\in M$ is said to be a \emph{one sided quasi-normalizer} of $B$ if there exists a finite set of elements $\{x_1,\ldots,x_n\}\subseteq M$ such that $Bx\subseteq \sum_{i=1}^nx_iB$. The set of one sided quasi-normalizers of $B$ in $M$ is denoted by $q\N^{(1)}_M(B)$ while the von Neumann algebra it generates is written $W^*(q\N^{(1)}_M(B))$ and  called the one sided quasi-normalizer algebra of $B$. We now present some immediate consequences of Theorem~\ref{T:main theorem}.

\begin{Corollary}\label{C:one sided quasi normalizers}
 The triple $B\subseteq W^*(q\N^{(1)}_M(B))\subseteq M$ has the relative weak asymptotic homomorphism property. 
\end{Corollary}

For the next corollary, we note that $B\subseteq M$ has the weak asymptotic homomorphism property precisely when the triple $B\subseteq B\subseteq M$ has the relative version.

\begin{Corollary}\label{C:WAHP}
 A von Neumann subalgebra $B$ of a finite von Neumann algebra $M$ has the weak asymptotic homomorphism property if and only if $W^*(q\N^{(1)}_M(B))=B$.
\end{Corollary}

Suppose that $B$ is a subfactor of a factor $M$ and $[M:N]<\infty$. Then $M=W^*(q\N^{(1)}_M(B))$ by  \cite[Proposition 1.3]{P-P}.
Thus we have the following corollary, which was first proved by
Grossman and Wiggins~\cite{G-W}.

\begin{Corollary}
 If $B$ is a finite index subfactor of a type ${\rm II}_1$ factor $M$ and $B\neq M$, then $B$ does not have the weak asymptotic homomorphism property.
\end{Corollary}

In comparing  $W^*(q\N^{(1)}_M(B))$ with the von Neumann algebra $q\N_M(B)''$ generated by the set of quasi-normalizers, it is clear that $W^*(q\N^{(1)}_M(B))\supseteq q\N_M(B)''$. It is an interesting question to know under what conditions  equality holds.
In this direction, we have the following result.
\begin{Proposition}
 If $B$ is an atomic von Neumann subalgebra of $M$, then 
\begin{equation}
W^*(q\N^{(1)}_M(B))=q\N_M(B)''=M.
\end{equation}
\end{Proposition}
\begin{proof}
 We need only show that $q\N_M(B)''=M$.  Since $B$ is atomic, $B=\oplus_{n=1}^NB_n$, where each $B_n$ is a full matrix algebra and $1\leq N\leq \infty$. Let $p_n$ be the central projections in $B$ corresponding to $B_n$. In the following we will show that $p_nMp_m\subseteq q\N(B)$ for $n\neq m$, which implies that $q\N_M(B)''=M$. Let $x\in p_nMp_m$. With respect to a choice of  matrix units of $B_n=p_nBp_n\cong M_r(\cc)$ and $B_m=p_mBp_m\cong M_s(\cc)$, we can write $x=(x_{ij})_{1\leq i\leq r,1\leq j\leq s}$. Let $y_{ij}$ be the $r\times s$ matrix with the $(i,j)$-th entry $x_{ij}$ and other entries 0 with respect to the same matrix units of $B_n$ and $B_m$. Now
\begin{align}
 Bx&=B_nx=\{(\lambda_{ij}x_{ij})_{1\leq i\leq r,1\leq j\leq s}: \lambda_{ij}\in\cc\}\notag\\
&=\sum_{1\leq i\leq r,1\leq j\leq s}y_{ij}B_m=\sum_{1\leq i\leq r,1\leq j\leq s}y_{ij}B.
\end{align}
By symmetry, $xB\subseteq \sum_{i=1}^n Bx_i$ for a finite set of elements $\{x_1,\ldots,x_n\}\subseteq M$. Thus $x\in q\N_M(B)$, completing the proof.
\end{proof}

Using Chifan's theorem in~\cite{Ch}, we have the following corollary of Theorem \ref{T:main theorem}. Note that the equality of the first and third algebras is already known by measure theoretic methods \cite{P-S}.
\begin{Corollary}\label{C:chifan}
 If $B$ is a masa in a separable type ${\rm II}_1$ factor $M$, then 
\begin{equation}
\N_M(B)''=W^*(q\N^{(1)}_M(B))=q\N_M(B)''.
\end{equation}
\end{Corollary}

In reference to Corollary \ref{C:chifan}, we do not know if the stronger equality 
$q\N^{(1)}_M(B)=q\N_M(B)$ holds for masas, even in the special cases  considered in Section 5.

We end this section with the following observation.

\begin{Theorem}\label{T:one=two}
 Let $N=W^*(q\N^{(1)}_M(B))$. If $p\in B'\cap \langle M,e_B\rangle$ is a finite projection in $\langle M,e_B\rangle$, then $p\leq e_N$. Furthermore, $W^*(q\N^{(1)}_M(B))=q\N_M(B)''$ if and only if $e_N$ is the supremum of all projections $p\in B'\cap \langle M,e_B\rangle$ such that $p$ is finite in $\langle M,e_B\rangle$.
\end{Theorem}
\begin{proof}
 The first statement is implied by the proof of Theorem~\ref{T:main theorem}. Suppose that $e_N$ is the supremum of all projections $p\in B'\cap \langle M,e_B\rangle$ such that $p$ is finite in $\langle M,e_B\rangle$. Then $e_N\left(B'\cap \langle M,e_B\rangle\right)e_N$ is a semi-finite von Neumann algebra.  Let $Q=q\N_M(B)''$. Clearly, $e_Q\leq e_N$, so suppose that $e_Q\neq e_N$. Then there is a nonzero finite projection $p\leq e_N-e_Q$ such that $p\in B'\cap \langle M,e_B\rangle$ and $p$ is finite in $\langle M,e_B\rangle$. By Lemma 1.4.2 of~\cite{Po2}, any projection $p'\in B'\cap\langle M,e_B\rangle$ with $p'\leq JpJ$ must be infinite. On the other hand, $JpJ\leq Je_NJ=e_N$ and therefore $JpJ\left(B'\cap \langle M,e_N\rangle\right) JpJ$ is semifinite. This is a contradiction. If $W^*(q\N^{(1)}_M(B))=q\N_M(B)''$, then $e_N$ is the supremum of all projections $p\in B'\cap \langle M,e_B\rangle$ such that $p$ is finite in $\langle M,e_B\rangle$ by 
\cite[Lemma 1.4.2 (iii)]{Po2}.
\end{proof}

\section{Group von Neumann algebras}\label{sec5}

In this section we will apply our previous results to the study of inclusions $L(H)\subseteq L(G)$ arising from inclusions $H\subseteq G$ of discrete groups. We will make the standard abuse of notation and write $g$ for a unitary in $L(G)$ and for a vector in $\ell^2(G)$. Thus we denote the Fourier series of $x\in L(G)$ by $x=\sum_{g\in G} \alpha_g g$ where $\sum_{g\in G} |\alpha_g|^2 < \infty$. We do not assume that $G$ is I.C.C., so that $L(G)$ may not be a factor. However, when using a trace, it will always be the standard one given by $\tau(e)=1$ and $\tau(g)=0$ for $g\in G\setminus \{e\}$.

The notion of one sided quasi-normalizers of von Neumann algebras has an obvious counterpart for group inclusions $H\subseteq G$. We say that $g\in G$ is a {\emph{one sided quasi-normalizer}} of $H$ if there exists a finite set $\{g_1,\ldots,g_n\}\subseteq G$ such that
\begin{equation}\label{eq5.1}
Hg\subseteq \cup_{i=1}^n \,g_iH.
\end{equation}
It is immediate that these elements form a semigroup inside $G$, denoted $q\mathcal{N}_G^{(1)}(H)$. However, for ease of notation, we will also denote this by $\Gamma$ throughout the section. There are two distinguished subgroups of $G$ associated with $\Gamma$. We denote by $H_1$ the maximal subgroup $\Gamma\cap \Gamma^{-1}$ inside $\Gamma$ (corresponding to the quasi-normalizers $q\mathcal{N}_G(H)$ defined by a two sided version of \eqref{eq5.1}). We let $H_2$ denote the subgroup of $G$ generated by $\Gamma$, and we note that the containment $H_1\subseteq H_2$ can be strict, as we show by a subsequent example. Many of the results in this section will depend on the following.
\begin{Theorem}\label{pro5.1}
Let $H\subseteq G$ be an inclusion of discrete groups, let $x\in q\mathcal{N}_{L(G)}^{(1)}(L(H))$, and write $x=\sum_{g\in G}\alpha_g g$ for its Fourier series. If $g_0\in G$ is such that $\alpha_{g_0}\ne 0$, then $g_0\in \Gamma$.
\end{Theorem}
\begin{proof}
Let $M=L(G)$ and $B=L(H)$. We may assume that $\|x\|=1$ and $Bx\subseteq \sum_{i=1}^r x_iB$ for a finite number of elements $x_1,\ldots,x_r\in M$. Let $\H$ be the closure of $BxB\xi$ in $L^2(M)$ so that  $\H$ is a $B$-bimodule. Since $\H\subseteq L^2(\sum_{i=1}^r x_iB)$, $\H$ is a finitely generated right $B$-module,  so there exist vectors $\eta_1,\ldots,\eta_k\in \H\subseteq L^2(M)$ such that 
\begin{equation}
 \eta=\sum_{i=1}^k\eta_i{\mathbb{E}}_B(\eta_i^*\eta),\quad   \eta\in\H,
\end{equation}
where $\eta,\eta_i$ are viewed as unbounded operators affiliated with $M$. In particular, we have
\begin{equation}\label{E:roger 1}
 bx=\sum_{i=1}^k \eta_i{\mathbb{E}}_B(\eta_i^*bx),\quad  b\in B.
\end{equation}

Set $C=\max\,\{\|\eta_i\|_2:1\leq i\leq k\}$, and let  $\eta_i=\sum_{g\in G}\alpha_g^i g$ be the Fourier series for $\eta_i$, $1\leq i\leq k$.  Since 
\begin{equation}
 \sum_{i=1}^k\sum_{g\in G}|\alpha_g^i|^2=\sum_{i=1}^k\|\eta_i\|_2^2<\infty,
\end{equation}
there is a finite set $S=\{g_1,\ldots, g_n\}\subseteq G$ such that
\begin{equation}\label{E:roger 2}
 k^2C^2\left(\sum_{i=1}^k\sum_{g\in S^c}|\alpha_g^i|^2\right)<|\alpha_{g_0}|^2.
\end{equation}

For each $h\in H$, it follows from~\eqref{E:roger 1} that
\begin{equation}
 L_{h}x=\sum_{i=1}^k \eta_i{\mathbb{E}}_B(\eta_i^*L_hx).
\end{equation}
Since ${\mathbb{E}}_B(\eta_i^*L_hx)\in L^2(B)$ for $1\leq i\leq k$, these elements have Fourier series which we 
write as ${\mathbb{E}}_B(\eta_i^*L_hx)=\sum_{h'\in H}\beta_{h'}^i h'$. Then
\begin{equation}\label{eq5.pq}
 \sum_{g\in G} \alpha_g hg=\sum_{i=1}^k\left(\sum_{g'\in G}\alpha_{g'}^ig'\sum_{h'\in H}\beta_{h'}^i h'\right).
\end{equation}
Comparing the coefficients of $hg_0$ on both sides of \eqref{eq5.pq}, we have 
\begin{equation}
 \alpha_{g_0}=\sum_{i=1}^k\left(\sum_{h'\in H}\alpha_{hg_0(h')^{-1}}^i\beta_{h'}^i\right).
\end{equation}
Since $\|x\|=1$, 
\begin{equation}
\|{\mathbb{E}}_B(\eta_i^*L_hx)\|_2\leq \|\eta_i^*L_hx\|_2\leq \|\eta_i\|_2\leq C,
\end{equation}
 and so $\sum_{h'\in H}|\beta_{h'}^i|^2\leq C^2$ for $1\leq i\leq k$. The Cauchy--Schwarz inequality gives
\begin{align}
 |\alpha_{g_0}|^2&\leq k^2 \sum_{i=1}^k\left(\sum_{h'\in H}\alpha_{hg_0(h')^{-1}}^i\beta_{h'}^i\right)^2\notag\\
&\leq k^2 \sum_{i=1}^k\left(\sum_{h'\in H}|\alpha_{hg_0(h')^{-1}}^i|^2\sum_{h'\in H}|\beta_{h'}^i|^2\right)\notag\\
&\leq k^2C^2 \sum_{i=1}^k\left(\sum_{h'\in H}|\alpha_{hg_0(h')^{-1}}^i|^2\right).
\end{align}
If $hg_0(h')^{-1}\in S^c$ for all $h'\in H$, then we have 
\begin{equation}
 |\alpha_{g_0}|^2\leq k^2C^2 \left(\sum_{i=1}^k\sum_{g\in S^c}|\alpha_g^i|^2\right),
\end{equation}
and this contradicts~\eqref{E:roger 2}. Thus there exists an $h'\in H$ such that $hg_0(h')^{-1}\in \{g_1,\ldots,g_n\}$, from which it follows that $hg_0\in g_iH$ for some $i$, $1\leq i\leq n$. Since $h\in H$ is arbitrary, we have shown that $Hg_0\subseteq \cup_{i=1}^n \,g_iH$, and therefore $g_0\in \Gamma$.
\end{proof}

A consequence of Theorem \ref{pro5.1} is that we can now describe both $q\mathcal{N}_{L(G)}(L(H))''$ and $W^*(q\mathcal{N}_{L(G)}^{(1)}(L(H)))$ in terms of groups.
\begin{Corollary}\label{cor5.4}
Let $H\subseteq G$ be an inclusion of discrete groups. Then
\begin{itemize}
\item[\rm{(i)}]
$q\mathcal{N}_{L(G)}(L(H))''=L(H_1)$;
\item[\rm{(ii)}]
$W^*(q\mathcal{N}_{L(G)}^{(1)}(L(H)))=L(H_2)$.
\end{itemize}
\end{Corollary}
\begin{proof}
The inclusion ``$\,\supseteq\,$'' is obvious in both cases. If $x\in q\mathcal{N}_{L(G)}^{(1)}(L(H))$ with Fourier series $\sum_{g\in G}\alpha_g g$, then any $g\in G$ for which $\alpha_g\ne 0$ must lie in $\Gamma\subseteq H_2$ by Theorem \ref{pro5.1}. This establishes ``$\,\subseteq\,$'' in (ii).

Now assume that $x\in q\mathcal{N}_{L(G)}(L(H))$, which is equivalent to $x,x^*\in q\mathcal{N}_{L(G)}^{(1)}(L(H))$. If $x=\sum_{g\in G}\alpha_g g$ then $x^*=\sum_{g\in G}\overline{\alpha_g}g^{-1}$, so Theorem \ref{pro5.1} gives $g,g^{-1}\in \Gamma$ whenever $\alpha_g\ne 0$. Then such elements $g$ lie in $\Gamma\cap\Gamma^{-1}=H_1$ and this shows the containment ``$\,\subseteq\,$'' in (i).
\end{proof}

Based on the above corollary, we can now present an example where the quasi-normalizers and the one sided quasi-normalizers are distinct.

\begin{Example}\label{ex5.3}
 \emph{Consider the free group ${\mathbb{F}}_\infty$, where the generators are written $\{g_i:\,i\in \mathbb{Z}\}$, and for each $n\in \mathbb{Z}$, let $K_n$ be the subgroup generated by $\{g_i:\,i\geq n\}$. The shift $i\rightarrow i+1$ on $\mathbb{Z}$ induces an automorphism $\phi$ of ${\mathbb{F}}_\infty$ defined on generators by $\phi(g_i)=g_{i+1}$, $i\in\mathbb{Z}$. Then $n\rightarrow \phi^n$ gives a homomorphism $\alpha:\mathbb{Z}\rightarrow {\rm Aut}({\mathbb{F}}_\infty)$, and we let $G$ be the semidirect product ${\mathbb{F}}_\infty\rtimes_\alpha\mathbb{Z}$.  Let $H=K_0$. In the following we will show that $H_1\neq H_2$. We denote by $t$ the generator of $\mathbb{Z}$.  Then every element of $G$ can be written as $w t^n$, where $w\in {\mathbb{F}}_\infty$. Note that $t H t^{-1}=K_1\subseteq H$. So $ Ht^{-1}\subseteq t^{-1} H$ and $t^{-1}\in H_2$. Suppose that $t^{-1}\in H_1$. Then $t^{-1}H\subseteq \cup_{i=-N}^N Ha_it^i$ for some large positive integer $N$ and some $a_i\in {\mathbb{F}}_\infty$. Multiplying on the right by $t$ gives $K_{-1}\subseteq
\cup_{i=-N}^N\,K_0a_it^{i+1}$ and so $K_{-1}\subseteq K_0a_{-1}$. If $r$ is
the total number of occurrences of $g_{-1}$ in $a_{-1}$, then
$g_{-1}^{r+1}\in K_{-1}$ but $g_{-1}^{r+1}\notin K_0a_{-1}$ and we reach a
contradiction. Thus $t^{-1}\notin H_1$ and so $H_1\ne H_2$.$\hfill\square$
}
\end{Example}

We now list some algebraic conditions on group inclusions $H\subseteq G$ that will be useful subsequently. The first two come from \cite{Di}. When $H$ is abelian, (C1) below gives a necessary and sufficient condition for $L(H)$ to be a masa in $L(G)$, while (C1) and (C2) combined give a sufficient condition for $L(H)$ to be a singular masa \cite{Di}. Subsequently (C2) alone was shown to be a necessary and sufficient condition in \cite{J-S} (see also the review of this paper,  MR2465603 (2010b:46127), by Stuart White).
\begin{itemize}
\item[{}]
\begin{itemize}
\item[\hspace{.25in}(C1)]
For each $g\in G\setminus H$, $\{hgh^{-1}:h\in H\}$ is infinite.
\item[\hspace{.25in}(C2)]
Given $g_1,\ldots,g_n \in G\setminus H$, there exists $h\in H$ such that
\begin{equation*}
g_ihg_j \notin H,\ \ \ 1\leq i,j\leq n.
\end{equation*}
\item[\hspace{.25in}(C3)]
If $g\in G$ and there exists a finite set $\{g_1,\ldots,g_n\}\subseteq G$ such that
\begin{equation*}
Hg\subseteq \cup_{i=1}^n \,g_iH,
\end{equation*}
then $g\in H$. ($\Gamma=q\mathcal{N}_G^{(1)}(H)=H$ in our notation).
\end{itemize}
\end{itemize}

We note that (C1) is a consequence of (C2) and also of (C3): if an element $g\in G\setminus H$ had only a finite number of $H$-conjugates $\{g_1,\ldots,g_n\}$, then (C2) would fail for the finite set
$\{g_1^{\pm 1},\ldots,g_n^{\pm 1}\}$, while (C3) would fail since we would have $Hg\subseteq \cup_{i=1}^{n}\,g_iH$. No abelian hypothesis on $H$ is required for this.

Combining Theorem~\ref{T:main theorem} and Corollary~\ref{cor5.4}, we obtain a purely algebraic characterization for the weak asymptotic homomorphism property. Note that we are not assuming $H$ to be abelian.
\begin{Corollary}\label{cor5.2}
Let $H\subseteq G$ be an inclusion of discrete groups. Then $L(H)\subseteq L(G)$ has the weak asymptotic homomorphism property if and only if condition {\rm{(C3)}} is satisfied.
\end{Corollary}

As mentioned above, condition (C2) is necessary and sufficient   to imply that $L(H)$ is a singular masa in $L(G)$ when $H$ is abelian (see~\cite{RSS,J-S}).  The following gives a different necessary and sufficient condition for singularity of $L(H)\subseteq L(G)$ in terms of the group structure. After Corollary \ref{cor5.3} has been proved, it will be apparent that conditions (C2) and (C3) are equivalent when $H$ is abelian. The direction (C2) $\Rightarrow$ (C3) is routine, but we do not have a purely group theoretic argument for the reverse implication.
\begin{Corollary}\label{cor5.3}
Let $H\subseteq G$ be an inclusion of discrete groups with $H$ abelian. Then $L(H)$ is a singular masa in $L(G)$ if and only if condition {\rm{(C3)}} is satisfied.
\end{Corollary}
\begin{proof}
Suppose that $L(H)$ is a singular masa in $L(G)$. From \cite{SSWW}, the inclusion $L(H)\subseteq L(G)$ has the weak asymptotic homomorphism property, so it is immediate from the definition that the triple $L(H)\subseteq L(H)\subseteq L(G)$ has the relative form. Theorem \ref{T:main theorem} then gives $q\mathcal{N}_{L(G)}^{(1)}(L(H))\subseteq L(H)$, and so condition (C3) holds.

Conversely, suppose that condition (C3) is valid.  Then  $q\mathcal{N}_{L(G)}^{(1)}(L(H))\subseteq L(H)$ follows from Theorem \ref{pro5.1}, and the weak asymptotic homomorphism property holds for $L(H)\subseteq L(G)$ by Theorem \ref{T:main theorem}. As noted before Corollary \ref{cor5.2}, condition (C1) is a consequence of condition (C3), so $L(H)$ is a masa in $L(G)$. Singularity now follows from \cite{SSWW}. 
\end{proof}

In the case that $L(H)$ is a masa in $L(G)$, we can now describe $\mathcal{N}_{L(G)}(L(H))''$ in terms of the normalizer $\mathcal{N}_G(H):=\{g\in G:gHg^{-1}=H\}$ at the group level. For this we need a preliminary group theoretic result.
\begin{Lemma}\label{lem5.5}
Let $H\subseteq G$ be an inclusion of discrete groups with $H$ abelian, and suppose that condition {\rm{(C1)}} holds. Let $g\in G$ be such that there exists a finite set $\{g_1,\ldots,g_n\}\subseteq G$ satisfying
\begin{equation}\label{eq5.b}
Hg\subseteq \cup_{i=1}^n \,g_iH.
\end{equation}
Then  $g\in \mathcal{N}_G(H)$.
\end{Lemma}
\begin{proof}
We may assume that the left cosets in \eqref{eq5.b} are a minimal set for which \eqref{eq5.b} holds. Thus they are pairwise distinct, and so disjoint, and minimality implies that
\begin{equation}\label{eq5.c}
Hg\cap g_iH\ne \emptyset,\ \ \ 1\leq i\leq n.
\end{equation}
Since there exist $h\in H$ and some integer $i$ such that $g=g_ih$, we may replace $g_i$ by $g_ih$ and renumber to further assume that $g=g_1$.

Now \eqref{eq5.b} implies that $HgH\subseteq \cup_{i=1}^n \,g_iH$ while the reverse containment follows from \eqref{eq5.c}. Since $HgH$ is invariant under left multiplication by elements $h\in H$, we obtain a representation $\pi$ of $H$ into the permutation group of $\{1,\ldots,n\}$ by defining $\pi_h(i)$ to be that (unique) integer $j$ so that $hg_i\in g_jH$, $1\leq i\leq n$. Let $K\subseteq H$ be the kernel of $\pi$, a finite index subgroup of $H$. Since $g=g_1$, we see that $kg\in gH$ for all $k\in K$. Let $\alpha \in {\rm{Aut}}\,(G)$ be defined by $\alpha(r)=g^{-1}rg$, for $r\in G$. Then, by definition of $K$, we have $\alpha(K)\subseteq H$. Thus
\begin{equation}\label{eq5.d}
\alpha(K)\subseteq H\cap\alpha(H)\subseteq \alpha(H),
\end{equation}
and so $K_1:=H\cap\alpha(H)$ has finite index in $\alpha(H)$ and satisfies
\begin{equation}\label{eq5.e}
K\subseteq \alpha^{-1}(K_1)\subseteq H.
\end{equation}
Then $\alpha^{-1}(K_1)$ has finite index in $H$, so we may list the cosets as $\alpha^{-1}(K_1)h_1,\ldots,\alpha^{-1}(K_1)h_m$ for some integer $m$ and elements $h_1,\ldots,h_m\in H$. For any $h\in H$ and $k\in K_1$,
\begin{equation}\label{eq5.f}
\alpha^{-1}(k)h_i\alpha^{-1}(h)h_i^{-1}\alpha^{-1}(k^{-1})
=h_i\alpha^{-1}(h)h_i^{-1},\ \ \ 1\leq i\leq m,
\end{equation}
since $\alpha^{-1}(k)$ commutes with both $h_i$ and $\alpha^{-1}(h)$. Thus $\alpha^{-1}(h)$ has only a finite number of $H$-conjugates, showing that $\alpha^{-1}(h)\in H$ from the hypothesis that condition (C1) holds. Thus $gHg^{-1}\subseteq H$.  But condition (C1) implies that $H$ (and hence $gHg^{-1}$) is maximal abelian in $G$, showing that $gHg^{-1}=H$. It follows that $g\in \mathcal{N}_G(H)$.
\end{proof}

\begin{Corollary}\label{cor5.6}
Let $H\subseteq G$ be an inclusion of discrete groups with $H$ abelian, and satisfying condition {\rm{(C1)}}, so that $L(H)$ is a masa in $L(G)$. Then
\begin{equation}\label{eq5.w}
\mathcal{N}_{L(G)}(L(H))''=L(\mathcal{N}_G(H)).
\end{equation}
In particular, $L(H)$ is a singular masa if and only if $\mathcal{N}_G(H)=H$, and is Cartan precisely when $H$ is a normal subgroup of $G$.
\end{Corollary}
\begin{proof}
By Theorem \ref{pro5.1}, any $u\in \mathcal{N}_{L(G)}(L(H))$ lies in $L(\Gamma)$. By Lemma \ref{lem5.5}, $\Gamma\subseteq \mathcal{N}_G(H)$, so $\mathcal{N}_{L(G)}(L(H))''\subseteq L(\mathcal{N}_G(H))$. Since the reverse inclusion is true for any subgroup $H$, the result follows.
\end{proof}

In specific cases this corollary is easy to apply. The properties of the various types of masas presented in \cite{Di} or in \cite[\S{2.2}]{S-S2} can now be verified trivially by using the equality of \eqref{eq5.w}. We also note that Corollary \ref{cor5.6} solves a question posed in~\cite[Remark 5.5]{SWW}.

\begin{Remark}\label{rem5.8}
\emph{The results of this section can be extended to the more general setting of inclusions $N\rtimes_{\theta}H\subseteq N\rtimes_{\theta}G$ where $N$ is a finite von Neumann algebra with a faithful normal trace $\tau$, $H\subseteq G$ are discrete groups, and $\theta$ is an action of $G$ on $N$ by trace preserving automorphisms. We make no assumptions that $G$ acts either freely or ergodically. The analog of Corollary \ref{cor5.4} is then the relations
\begin{equation}\label{eq5.g}
q\mathcal{N}_{N\rtimes_{\theta}G}(N\rtimes_{\theta}H)''=N\rtimes_{\theta}H_1,
\ \ W^*(q\mathcal{N}_{N\rtimes_{\theta}G}^{(1)}(N\rtimes_{\theta}H))=
N\rtimes_{\theta}H_2,
\end{equation}
 which are seen to be generalizations by taking $N=\mathbb{C}1$ and $\theta$ the trivial action. We omit the details since they are so similar to what has already been presented, and we mention only the one small change that is necessary. The Fourier series $\sum_{g\in G} \alpha_g g$ of Theorem \ref{pro5.1} is replaced by $\sum_{g\in G}x_g g$ with $x_g\in N$ and $\sum_{g\in G}\|x_g\|_2^2<\infty$, and $\|x_g\|_2$ is substituted in all calculations involving
$|\alpha_g|$.$\hfill\square$}
\end{Remark}

\section{Concluding remarks}
Let $B$ be a von Neumann subalgebra of $M$. Various notions of ``$\,$normalizers$\,$'' have been introduced:
\begin{itemize}
 \item[(i)] normalizers $\N_M(B)$ (\cite{Di}): a unitary operator $u\in M$ is a normalizer of $B$ if $uBu^*=B$;
 \item[(ii)] one sided normalizers ${\mathcal {ON}}_M(B)$ (\cite{SWW}): a unitary operator $u\in M$ is a one sided normalizer of $B$ if $uBu^*\subseteq B$;
\item[(iii)] groupoid normalizers ${\mathcal {GN}}_M(B)$ (\cite{Dye}):  a partial isometry $v\in M$ is a groupoid normalizer of $B$ if $vBv^*\subseteq B$ and $v^*Bv\subseteq B$;
\item[(iv)] intertwiners ${\mathcal {GN}}_M^{(1)}(B)$ (\cite{FSWW}):  a partial isometry $v\in M$ is an intertwiner of $B$ if $v^*v\in B$ and $vBv^*\subseteq B$;
\item[(v)] quasi-normalizers $q{\mathcal N}_M(B)$ (\cite{Po1.5}):  an operator $x\in M$ is a quasi-normalizer of $B$ if there exists a finite number of elements $x_1,\ldots,x_n\in M$ such that $Bx\subseteq \sum_{i=1}^n x_iB$ and  $xB\subseteq \sum_{i=1}^n Bx_i$;
\item[(vi)] one sided quasi-normalizers $q\N^{(1)}_M(B)$: an operator $x\in M$ is a one sided quasi-normalizer of $B$ if there exist a finite number of elements $x_1,\ldots,x_n\in M$ such that $Bx\subseteq \sum_{i=1}^n x_iB$.
\end{itemize}
The relations between  von Neumann algebras generated by the above ``$\,$normalizers$\,$'' are the following:
\begin{equation}\label{E:diagram}
 \begin{array}{ccccc}
  \N_M(B)''&\subseteq & {\mathcal {GN}}_M(B)''&\subseteq & q{\mathcal N}_M(B)''\\
\begin{sideways}
 $\supseteq$
\end{sideways}
 &&\begin{sideways}
 $\supseteq$
\end{sideways}&&\begin{sideways}
 $\supseteq$
\end{sideways}\\
{\mathcal {ON}}_M(B)''&\subseteq& W^*({\mathcal {GN}}^{(1)}_M(B))&\subseteq& 
W^*(q\N^{(1)}_M(B))\\
 \end{array}.
\end{equation}

By Corollary \ref{C:chifan}, if $B$ is a masa in a type ${\rm II}_1$ factor $M$, then $\N_M(B)''=W^*(q\N^{(1)}_M(B))$ and therefore all of the above ``$\,$normalizer algebras$\,$'' are the same.  On the other hand, for each ``$\,X\subseteq Y\,$'' in the above diagram, there are examples of inclusions of finite von Neumann algebras such that $X\neq Y$ (see~\cite{FSWW,SWW}).  Among the above ``$\,$normalizer algebras$\,$'', $W^*(q\N^{(1)}_M(B))$ has the following two special properties. The first of these is the formula \eqref{eq6.pq} for tensor products. This is an outgrowth of the analogous formulas
\begin{equation}
{\mathcal {GN}}_{M_1\tensor M_2}\left(B_1\tensor B_2\right)''={\mathcal {GN}}_{M_1}(B_1)''
\tensor {\mathcal {GN}}_{M_2}(B_2)''
\end{equation}
and
\begin{equation}
W^*\left({\mathcal {GN}}^{(1)}_{M_1\tensor M_2}\left(B_1\tensor B_2\right)\right)=W^*\left({\mathcal {GN}}^{(1)}_{M_1}(B_1)\right)
\tensor W^*\left({\mathcal {GN}}^{(1)}_{M_2}(B_2)\right)
\end{equation}
established in \cite{FSWW} under the hypothesis that $B_i'\cap M_i\subseteq B_i$, and which can fail without some such assumption. In contrast, the next proposition requires no restrictions.
\begin{Proposition}\label{pro6.1}
 Let $B_i\subseteq M_i$ be  inclusions of finite von Neumann algebras, $i=1,2$. Then
\begin{equation}\label{eq6.pq}
 W^*\left(q\N^{(1)}_{M_1\tensor M_2}\left(B_1\tensor B_2\right)\right)=
W^*(q\N^{(1)}_{M_1}(B_1))\tensor W^*(q\N^{(1)}_{M_2}(B_2)).
\end{equation}
\end{Proposition}
\begin{proof}
 Suppose that $x_1\in M_1$ satisfies $B_1x_1\subseteq\sum_{i=1}^{n_1}y_iB_1$ for a finite number of elements $y_1,\ldots, y_{n_1}$ in $M_1$, and $x_2\in M_2$ satisfies $B_2x_2\subseteq\sum_{i=1}^{n_2}z_iB_2$ for a finite number of elements $z_1,\ldots, z_{n_2}$ in $M_2$.  Then $(B_1\tensor B_2)(x_1\otimes x_2)\subseteq \sum_{i=1}^{n_1}\sum_{j=1}^{n_2}(y_i\otimes z_j)(B_1\tensor B_2)$. This proves that
\begin{equation}
  W^*\left(q\N^{(1)}_{M_1\tensor M_2}\left(B_1\tensor B_2\right)\right)\supseteq W^*(q\N^{(1)}_{M_1}(B_1))\tensor W^*(q\N^{(1)}_{M_2}(B_2)).
\end{equation}
On the other hand, the triple $B_i\subseteq W^*(q\N^{(1)}_{M_i}(B_i))\subseteq M_i$ has the relative weak asymptotic homomorphism property by Corollary~\ref{C:one sided quasi normalizers}, $i=1,2$, and so 
\begin{equation}
B_1\tensor B_2\subseteq W^*(q\N^{(1)}_{M_1}(B_1))\tensor 
W^*(q\N^{(1)}_{M_2}(B_2))\subseteq M_1\tensor M_2
\end{equation}
 also has the relative weak asymptotic homomorphism property. By Theorem~\ref{T:main theorem}, 
we have the reverse containment
\begin{equation}
  W^*\left(q\N^{(1)}_{M_1\tensor M_2}\left(B_1\tensor B_2\right)\right)\subseteq W^*(q\N^{(1)}_{M_1}(B_1))\tensor W^*(q\N^{(1)}_{M_2}(B_2)),
\end{equation}
completing the proof.
\end{proof}

We end with  two results that discuss the situation of a cut down of $B\subseteq M$ to an inclusion $eBe\subseteq eMe$  for a projection $e\in B$.

\begin{Proposition}\label{P:reduction}
 Let $e\in B$ be a projection. Then $W^*(q\N^{(1)}_{eMe}(eBe))=eW^*(q\N^{(1)}_M(B))e$.
\end{Proposition}
\begin{proof}
We only prove that $e \left(q\N^{(1)}_M(B)\right) e\subseteq W^*(q\N^{(1)}_{eMe}(eBe))$. The proof of $q\N^{(1)}_{eMe}(eBe)\subseteq eW^*(q\N^{(1)}_M(B))e$ is similar. Suppose that $z$ is a central projection in $B$ such that $z=\sum_{j=1}^n v_jv_j^*$ with the $v_j$'s partial isometries in $B$ and $v_j^*v_j\leq e$. Write $e_0=ez$.  If $x\in M$ satisfies $Bx\subseteq \sum_{i=1}^r x_iB$, then
\begin{align}
eBe e_0xe_0&\subseteq eBzxe_0=ezBxe_0\subseteq e_0\sum_{i=1}^r x_iBe_0=e_0\sum_{i=1}^r x_izBe_0\notag\\
&=\sum_{i=1}^r\sum_{j=1}^n (e_0x_iv_j)(v_j^*Be_0)
\subseteq \sum_{i=1}^r\sum_{j=1}^n (e_0x_iv_j)(eBe).
\end{align}
Therefore, $e_0xe_0\in q\N^{(1)}_{eMe}(eBe)$. Since the central
support of $e$ in $B$ can be approximated arbitrarily well by such
special central projections $z$, $e_0$ approximates $e$ arbitrarily well,
 and $exe\in W^*(q\N^{(1)}_{eMe}(eBe))$.
\end{proof}

Combining Proposition~\ref{P:reduction} and Corollary~\ref{C:WAHP}, we obtain the following consequence.
\begin{Corollary}
 Suppose that $B$ has the weak asymptotic homomorphism property in $M$ and $e\in B$ is a projection. Then $eBe$ has the weak asymptotic homomorphism property in $eMe$.
\end{Corollary}

\vspace{.2in}

\noindent Junsheng Fang

\noindent Department of Mathematics, Texas A\&M University,

\noindent College Station, TX 77843

\noindent {\em E-mail address: } [Junsheng Fang]\,\,
jfang\@@math.tamu.edu

\vspace{.2in}

\noindent Mingchu Gao

\noindent Department of Mathematics, Louisiana College,

\noindent Pineville, LA 71359

\noindent {\em E-mail address: } [Mingchu Gao]\,\,
gao\@@lacollege.edu

\vspace{.2in}

\noindent Roger R. Smith

\noindent Department of Mathematics, Texas A\&M University,

\noindent College Station, TX 77843

\noindent {\em E-mail address: } [Roger R. Smith]\,\,
rsmith\@@math.tamu.edu

\end{document}